# The sample fraction in peaks-over-threshold problems where the second-order expansion is valid with specific reference to the generalized Pareto distribution

J. Martin van Zyl

*Department of Mathematical Statistics and Actuarial Science, University of the Free State, Bloemfontein, South Africa*

e-mail: wwjvz@ufs.ac.za

In samples from a heavy-tailed distribution the peaks-over-threshold approach is often applied and above a threshold the distribution of excesses over the threshold have asymptotically a generalized Pareto distribution. The tail function can be approximated by using a second-order expansion and based on the estimated parameters of the second–order approximation an expression for the optimal sample fraction of largest order statistics can be derived. Given that the observations are above the threshold, an expression is derived for the percentile above which the second-order approximation is valid.

## 1. Introduction

In samples from a heavy-tailed distribution the peaks-over-threshold approach (Smith and Davidson, 1990) is often applied and above a threshold the distribution of excesses over the threshold has asymptotically a generalized Pareto distribution (GPD), (Balkema and de Haan(1974), Pickands (1975)). Hall (1982), Hall and Welsh (1985) considered distribution functions with a positive index, $\alpha$, for which the tail function can be approximated by using a second-order expansion and based on the estimated parameters of the second–order approximation an expression for the optimal sample fraction and bias of the Hill estimator (Hill, 1975) can be derived. Many of the references on this topic can be found in the review paper by Gomes and Guillou (2015).

Let $\alpha$ denote the index of the excesses over the threshold. Given that the observations are above the threshold, it will be shown that the expansion is valid for observations above the $1-1/2^{\alpha}$ percentile. An estimate of $\alpha$ can easily be made by using a small number of the largest observations and applying the Hill estimator (Hill, 1975). This result can also be used as a lower bound for the threshold since it is the smallest value for which there is power behaviour in the tails of the distribution. The fact that the expansion is valid, does not mean it is a good when approximating the tail function, and to be a good approximation the order statistics should also be larger than one since the expansion is in terms of the inverse of the order statistics.

It is shown that the same condition hold for the Burr distribution.

The distribution function of the generalized Pareto distribution is

$F(x)=1-(1+(\xi/\sigma)x)^{-\alpha}$ and tail function

$$1-F(x)=(1+(\xi/\sigma)x)^{-\alpha},\ x>0, \qquad (1)$$

where $\xi = 1/\alpha$ the shape parameter, $\sigma$ a scale parameter, and $x$ an excess over a threshold $\mu$, or $x = z - \mu$. The second-order approximation introduced by Hall (1982), Hall and Welsh (1985) is written as

$$1 - F(x) = Ax^{-\alpha}\{1 + Bx^{-\beta} + o(x^{-\beta})\}, \; x \to \infty$$

$$\approx cx^{-a} + dx^{-b} \quad , \qquad (2)$$

where $c > 0, d \neq 0, b > a > 0$.

In section 2 it will be used to shown that the expansion is valid for a GPD, $F$, if

$$F(x) > 1 - 2^{-\alpha}. \qquad (3)$$

In practice $\alpha$ can be estimated using a small number of the largest excesses and the proportion calculated by using an approximate threshold. Figures are shown to illustrate the idea and in section 3 this will be applied to simulated samples for various values of $\alpha$.

For the Fréchet distribution it is trivial to show that the expansion is valid for all observations and any value of the index $\alpha$. For a standard Fréchet distribution it can be shown that by using the condition that observations must be larger than one, $P(X > 1)$, the observations above the upper $1 - 1/e$ percentile can be included for a good 2nd order approximation.

## 2. The second-order approximation applied to a GPD

A theorem will be derived to find the point where the approximation is valid in terms of percentiles.

Theorem: Let $Z_1,...,Z_n$ denote the largest observations of a sample and $X_1 = Z_1 - \mu,..., X_n - \mu$ denote the excesses above a threshold , $\mu$ . The threshold is large enough such that the asymptotic distribution of the excesses over the threshold have a GPD distribution $G$ with scale parameter $\sigma$ and index $\alpha$ . The second order approximation applied to the excess is valid for observations such that $F(x) > 1 - 1/2^{\alpha}$ . If $N$ observations are less than the threshold and $n - N$ larger or equal to the threshold, the expansion is valid for observations larger than the percentile $F(x) > (N/n) + ((n-N)/n)(1-1/2^{\alpha})$ .

Proof: Let $G(x) = 1 - F(x) = (x/\alpha\sigma)^{-\alpha}(1 + \alpha\sigma/x)^{-\alpha}$ . The binomial theorem can be applied to the term $(1 + \alpha\sigma/x)^{-\alpha}$ for $\alpha\sigma/x$ <1. Thus

$$\alpha\sigma/x < 1$$

$$x/\alpha\sigma > 1$$

$$1 + x/\alpha\sigma > 2$$

$$(1 + x/\alpha\sigma)^{\alpha} > 2^{\alpha},$$

Thus for

$$F(x) = 1 - (1 + x/\alpha\sigma)^{-\alpha} > 1 - 2^{-\alpha},$$

or $G(x) < 2^{-\alpha}$.

An easier way to derive this result is:

$$\begin{aligned}\mathrm{P}(X > \alpha\sigma) &= 1 - F(\alpha\sigma) \\ &= (1 + (1/\alpha\sigma)(\alpha\sigma))^{-\alpha} \\ &= 2^{-\alpha}.\end{aligned}$$

In terms of the observed sample it means that $Z_j \geq \mu + \alpha\sigma,\ j = 1,...n$.

The Burr distribution with

$$f(x) = \tau\alpha\lambda^{\alpha} x^{\tau-1} / (\lambda + x^{\tau})^{\alpha+1}$$

and $$\begin{aligned}1 - F(x) &= \lambda^{\alpha} / (\lambda + x^{\tau})^{\alpha} \\ &= (\lambda^{\alpha} / x^{-\alpha\tau})(1 + \lambda / x^{\tau})^{-\alpha}.\end{aligned}$$

The term $(1 + \lambda / x^{\tau})^{-\alpha}$ can be expanded into a series in terms of the inverse of $x$, and the expansion is valid for $\lambda^{1/\alpha} < x,\ F(x) > F(\lambda^{1/\tau}) = 1 - 2^{-\alpha}$. Thus the largest observations higher than this percentile must be used for the power expansion to be valid, and to be a good approximation for $x > 1$. The expansion with two terms is

$$1 - F(x) \approx \lambda^{\alpha} / x^{\alpha\tau} - \alpha\lambda^{\alpha} / x^{(\alpha+1)\tau}.$$

By making use of the binomial theorem the tail function can be expanded as series to infinity. The expansion is valid for $x > \alpha\sigma$. Denote the second order approximation by $G^*(x) \approx 1 - F(x)$.

$$
\begin{aligned}
1 - F(x) &= (x/\alpha\sigma)^{-\alpha}(1 + \alpha\sigma/x)^{-\alpha} \\
&= (x/\alpha\sigma)^{-\alpha}(1 + (-\alpha)(x/\alpha\sigma) + (1/2)(-\alpha)(-\alpha-1)(x/\alpha\sigma)^2 + \\
&\qquad (1/6)(-\alpha)(-\alpha-1)(-\alpha-2)(x/\alpha\sigma)^3 + \ldots \\
&\approx (\alpha\sigma)^{\alpha} x^{-\alpha} - \alpha^{(\alpha+2)}\sigma^{(\alpha+1)}x^{-(\alpha+1)}.
\end{aligned}
$$

The constants are $c = (\alpha\sigma)^{\alpha}$, $a = \alpha$, $d = -\alpha^{(\alpha+2)}\sigma^{(\alpha+1)}$, $b = \alpha + 1$, and $G^*(\alpha\sigma) = 1 - \alpha$. In the following figure the second order tail function approximation for a GPD with $\sigma = 1, \alpha = 2$, is shown as a function of $x$.

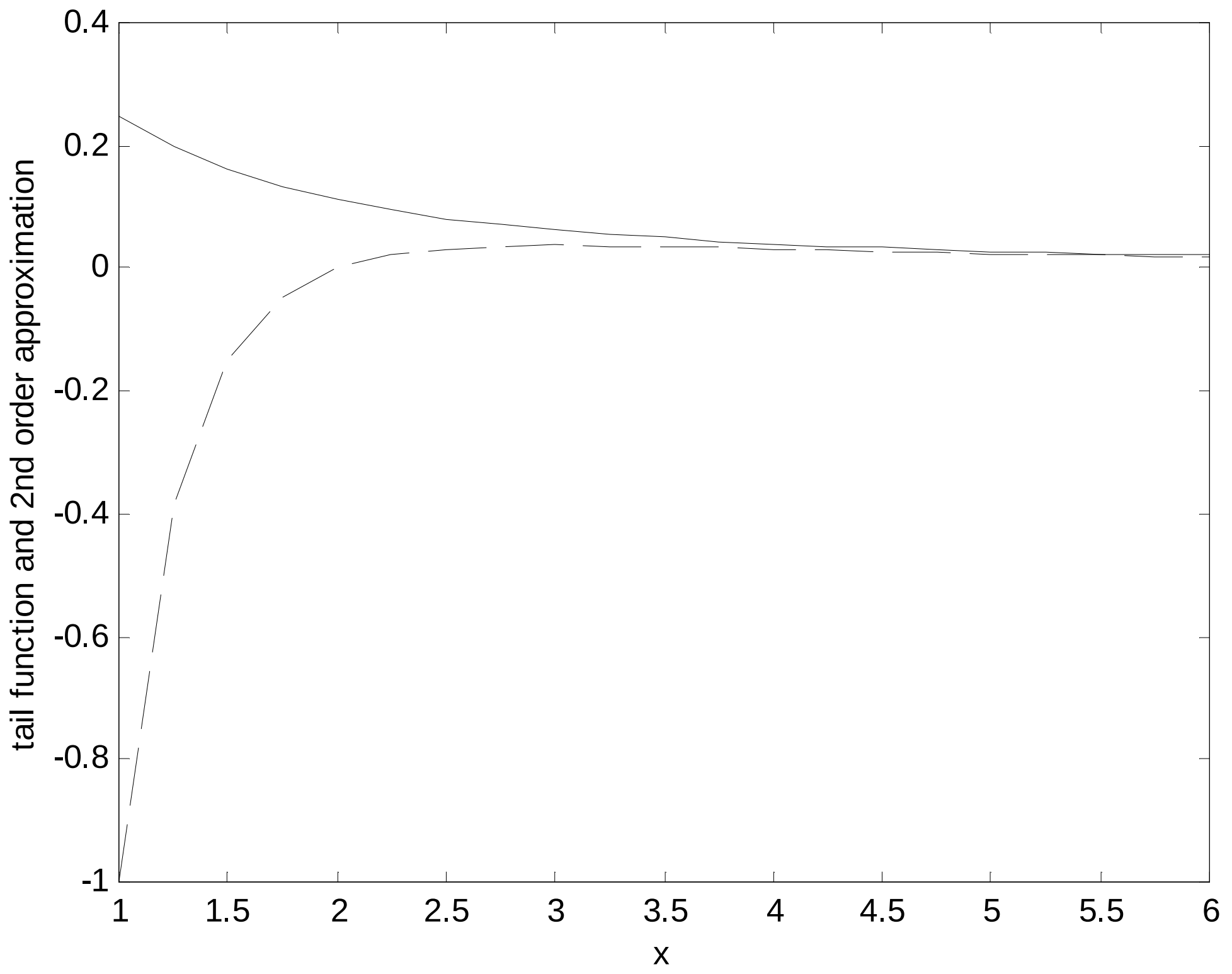

Figure 1. Second order approximation and true tail function a of GPD, $\alpha = 2, \sigma = 0.5$

In the following figure a sample of size $n = 250$ is simulated from a GPD with threshold zero, thus similar data as excesses over a threshold. The upper 25% of observations are used, say $m$, and $G^*(x_{(j)}), 1 - F(x_{(j)})$ plotted on the vertical axis against the empirical distribution function $(j-1)/m$ on the horizontal axis for $j = 1, ..., m$. The parameters of the GPD are $\sigma = 0.5, \alpha = 1$. As one would expect with $\alpha = 1$, the approximation should be valid for the upper $(1 - 2^{-\alpha})100 = 50\%$ of observations. Since the distribution is very heavy tailed these observations are much larger than one and one can see the second order approximation is not only valid but also a good approximation of the tail function.

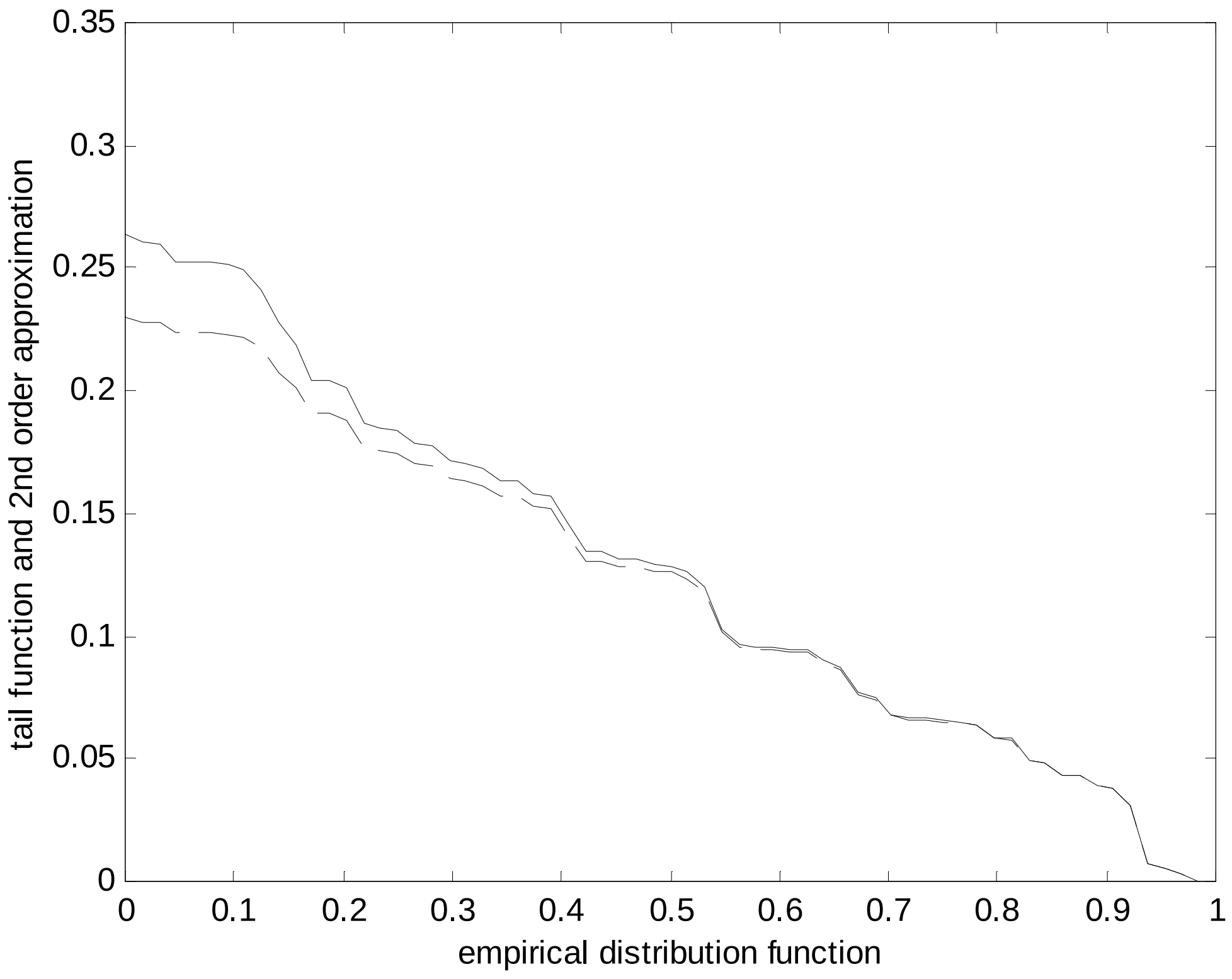

Figure 2. Approximate and true tail function of largest 25% of observations in a sample of with a GPD distribution with $\sigma = 0.5, \alpha = 1$ plotted against the empirical distribution function. The dashed line is the approximation and the solid line true tail function.

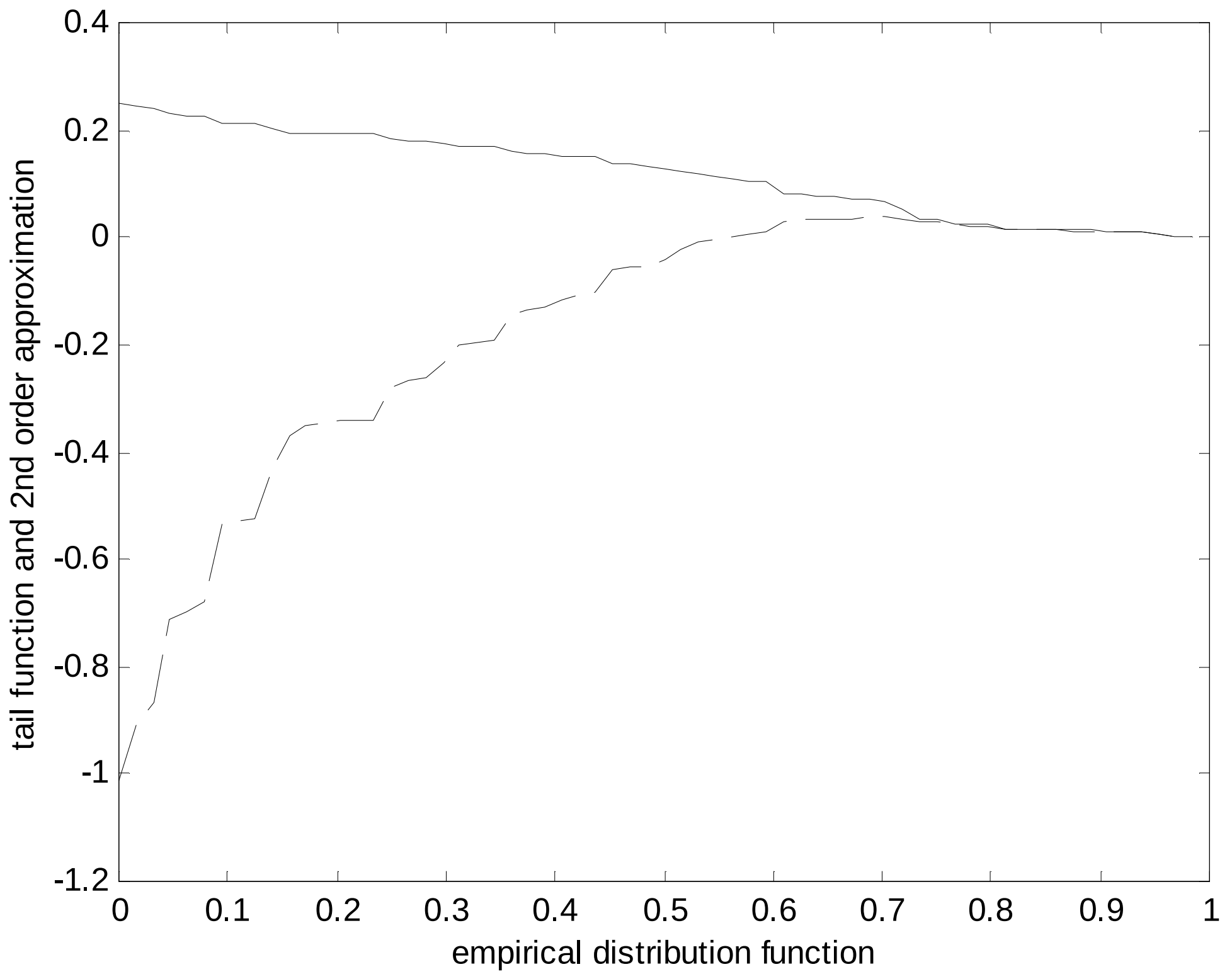

Figure 3. Approximate and true tail function of largest 25% of observations in a sample of with a GPD distribution with $\sigma = 2, \alpha = 2$ plotted against the empirical distribution function. The dashed line is the approximation and the solid line true tail function.

The Fréchet distribution with distribution function $F$ and index $\alpha$, $F(x) = e^{-x^{-\alpha}}$, the expansion is valid for any positive value of $\alpha$ and positive x, and it is

$$1 - F(x) \approx x^{-\alpha} - (1/2)x^{-2\alpha}.$$

It can be seen that this series would converge faster for $x > 1$ and especially for large values of $x$. The proportion of observations where $x > 1$ is $F(1) = e^{-1}$. A guideline would be to use the largest 63% of the observations.

For the $t$-distribution function with $\nu$ degrees of freedom and density function

$$f(t) = 1/(\sqrt{\nu} B(\nu/2, 1/2)(1+t^2/\nu)^{-(\nu+1)/2},$$

It follows that the tail function expansion is valid for $t > \sqrt{\nu}$

$$1 - F(x) = 1/(\sqrt{\nu} B(\nu/2, 1/2) \int_x^{\infty} (\nu/t^2)^{(\nu+1)/2} (1 - \frac{(\nu+1)}{2} \nu t^{-2} + ...) dt$$

$$\approx 1/(\sqrt{\nu} B(\nu/2, 1/2)(\nu^{(\nu+1)/2})(\nu^{-1} x^{-\nu} - 0.5\nu(\nu+1)/(\nu+2) x^{-(\nu+2)})$$

$$= 1/(B(\nu/2, 1/2)(\nu^{\nu/2-1} x^{-\nu} - 0.5\nu^{\nu/2+1}(\nu+1)/(\nu+2) x^{-(\nu+2)}) .$$

The tail function can be calculated exactly by making use of the transformation $x = \nu/(\nu + t^2)$ and the incomplete beta function. In the point $t = \sqrt{\nu}$ , $x = \nu/(\nu + (\sqrt{\nu})^2) = 0.5$, and

$$1 - F(\sqrt{\nu}) = 0.5(1/B(\nu/2, 0.5)) \int_0^{0.5} t^{\nu/2-1} t^{-0.5} dt .$$

This is tabulated in table 1 for some degrees of freedom. Notice that in practice if the negative values are taken into consideration the proportions would be twice those tabulated.

| Degrees of freedom $\nu$ | $P(X > \sqrt{\nu})$ |
|---|---|
| 1 | 0.25 |
| 2 | 0.1464 |
| 3 | 0.0908 |
| 4 | 0.0581 |
| 5 | 0.0378 |
| 6 | 0.0249 |
| 7 | 0.0166 |
| 8 | 0.0111 |
| 9 | 0.0075 |
| 10 | 0.0051 |

In the following figure a plot is made of the tail function of a t-distribution with degrees of freedom , $\nu = 2$, and a sample of size $m = 250$ simulated. The largest 25% of the observations are considered. Thus $\hat{F}_j = (j-1)/m$ against $1 - F(x_{(j)})$ and $G^*(x_{(j)})$, $j = 1,...,k$ for k observations considered. The solid line is the tail function and 2nd order approximation calculated in the points which were observed and plotted against the empirical distribution function.

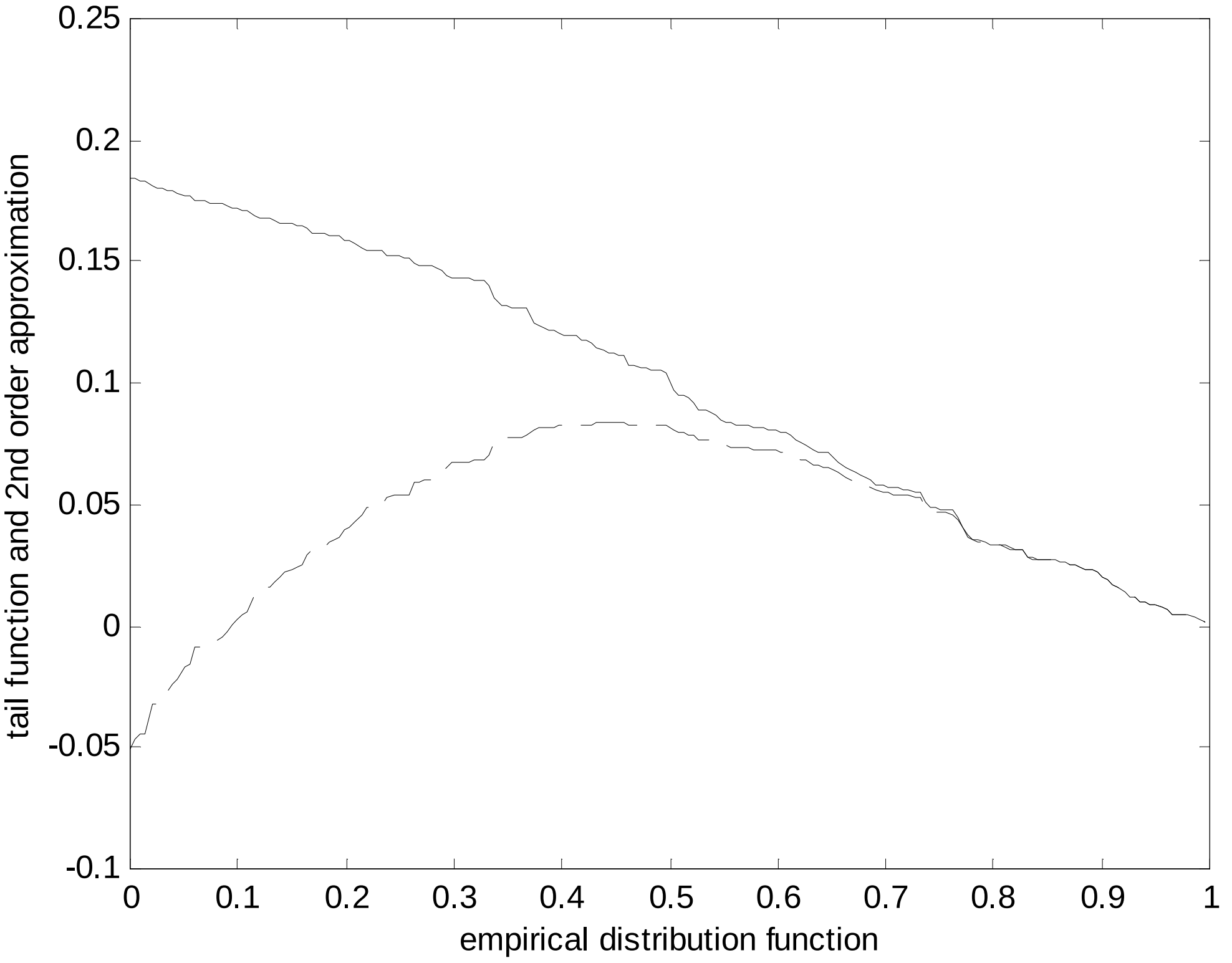

.

Figure 4. Approximate and true tail function of largest 20% of observations in a sample of size $m = 1000$ from a $t-$distribution with $\nu = 2$, plotted against the empirical distribution function. The dashed line is the approximation and the solid line true tail function.

## 3. Conclusions

A guideline to choose the fraction of observations above a threshold for which a second order approximation of the tail function is valid was derived. This can also be used as a conservative estimate of the threshold. DuMouchel (1983) suggested in the context of observations from a stable distribution, that not more than largest 10% observations should be included when using the Hill estimator. Some researchers, for example

Galbraith and Zernov (2004) use this rule when using the Hill estimator to measure the index of heavy-tailed data. In terms of the rule derived here 10% would correspond to $\alpha = 2.3$, and more than 10% of the largest observations can used for heavier tails.